\documentclass[3pt]{article}
\usepackage{graphics, amssymb, amsmath, mathrsfs,graphicx,theorem}
\newtheorem{theorem}{Theorem}[section]

\newtheorem{lemma}{Lemma}[section]

\newtheorem{definition}{Definition}[section]
\newtheorem{remark}{Remark}[section]

\newenvironment{proof}{{\noindent\it Proof:}\quad}{\hfill $\square$}
\numberwithin{equation}{section}

\usepackage[bottom]{footmisc}

\begin{document}

\title{Existence of solutions for first-order Hamiltonian stochastic impulsive differential equations with Dirichlet boundary conditions}

\author{Yu Guo$^{a}$\footnotemark[4] , Xiao-Bao Shu$^{a}$\thanks{ Corresponding author: Xiao-Bao Shu.
		Email: sxb0221@163.com(Xiao-Bao Shu).} 
	, Qian bao Yin$^{a}$\footnotemark[4]\\
\footnotesize{
	$^a$ College of Mathematics and Econometrics, Hunan University,}\\
\footnotesize{Changsha, Hunan 410082, PR China.}\\
}
\renewcommand{\thefootnote}{\fnsymbol{footnote}}
\footnotetext[4]{Yu Guo and Qian bao Yin contributed equally to this manuscript.}
\date{}
\maketitle
\begin{center}
	\begin{minipage}{139mm}
		
		\rule[4pt]{14.3cm}{0.05em}
		\noindent{\bf Abstract:}
		In this paper, we study the sufficient conditions for the existence of solutions of first-order Hamiltonian stochastic impulsive differential equations under Dirichlet boundary value conditions. By using the variational method, we first obtain the corresponding energy functional. And by using Legendre transformation, we obtain the conjugation of the functional. Then the existence of critical point is obtained by mountain pass lemma. Finally, we assert that the critical point of the energy functional is the mild solution of the first order Hamiltonian stochastic impulsive differential equation.Finally, an example
		are presented to illustrate the feasibility and effectiveness of our results.\\
		
		
		\rule[4pt]{14.3cm}{0.05em}
	\end{minipage}
\end{center}\vspace{5mm}

\section{Introduction \label {Sec1}}

\quad\quad The phenomenon of random impulse exists widely in nature, and differential equation is one of the most powerful mathematical tools in scientific research. The content of this paper is the modern nonlinear methods and applications of stochastic impulsive differential equations. The purpose of this study is to explore the existence of solutions for a class of stochastic impulsive differential equations by using the critical point theory and variational method in modern mathematics. Stochastic impulsive differential equations are used in many fields, computer, finance, biomedicine, artificial intelligence, optimal control model [19, 20], and so on. Therefore, the study of impulsive differential equations has greatly practical significance. But in some practical cases, such as the mechanical problem [11, 15, 23], the impulse is random, so that means that the solution of stochastic impulsive differential equation is a random process, which is different from the corresponding fixed impulsive differential equation whose solution is a piecewise continuous function. Many scholars have studied fixed impulsive differential equations [ 4, 5, 9, 10, 14, 18, 21, 22, 25], while stochastic impulsive differential equations [6, 8, 17, 27, 29] are rarely involved. It is also of great significance to explore the use of nonlinear methods in the study of differential equations. This paper lays a foundation for the application of these methods in various fields of differential equation research in the future, and fills the gap in the field of studying stochastic impulse differential equation by nonlinear method.

Many scholars focus on the existence and multiplicity of solutions, and get most results of existence of solutions[7, 12, 13, 26, 28, 30, 31]. 
For example, Ravi P. Agarwal proved the multiplicity of second order impulsive differential equations in [26] by using Leggett Williams fixed point theorem:
\begin{align*}
	\label{eq1}
	\left\{
	\begin{array}{ll}
		y''(t)+\phi(t)f(y(t))=0&for\quad t\in[0,1]\backslash\{t_1,\cdots,t_m\},\\
		\Delta y(t_k)=I_k(y(t_k^-)),&k=1,\cdots,m,\\
		\Delta y'(t_k)=J_k(y(t_k^-)),&k=1,\cdots,m,\\
		y(0)=y(1)=0.
	\end{array} 
	\right.
\end{align*}

Here $\Delta y(t_k)=y(t_k^+)-y(t_k^-)$ where $y(t_k^+)$ (respectively $y(t_k^-)$) denote the right limit (respectively left limit) of $y(t)$ at $t=t_k$. Also $\Delta y'(t_k)=y'(t_k^+)-y'(t_k^-)$.  
The upper and lower solution method is also used to study impulsive differential equations [30, 31]. In [30], Jianhua Shen and Weibing Wang established the existence condition of the solution by using the upper and lower solution method and Schauder's fixed point theorem:
\begin{align*}
	\left\{
	\begin{array}{ll}
		x''(t)=f(t,x(t),x'(t))& t\in J,t\neq t_k\\
		\Delta x(t_k)=I_k(x(t_k)),&k=1,2,\cdots,p,\\
		\Delta x'(t_k)=J_k(x(t_k),x'(t_k)),&k=1,2,\cdots,p,\\
	g(x(0),x'(0))=0,\quad h(x(1),x'(1))=0
	\end{array} 
	\right.
\end{align*}
where $J=[0,1],f:J\times R^2\to R$ is continuous, $I_k,J_k\in C(R)$ for $1\leq k\leq p,0=t_0<t_1<t_2<\cdots<t_p<t_{p+1}=1,\Delta x(t_k)=x(t_k^+)-x(t_k^-)$ denotes the jump of $x(t)$ at $t=t_k,x(t_k^+)$ and $x(t_k^-)$ represent the right and left limits of $x(t)$ at $t=t_k$ respectively, and $g,h:R^2\to R$ are continuous. $\Delta x'(t_k)=x'(t_k^+)-x'(t_k^-)$, where 
\begin{align*}
x'(t_k^-):=\lim\limits_{h\to 0^-}h^{-1}[x(t_k+h)-x(t_k)],\quad x'(t_k^+):=\lim\limits_{h\to 0^+}h^{-1}[x(t_k+h)-x(t_k)].
\end{align*}
Let $J^{\ast}=J\backslash\{t_1,t_2,\cdots,t_p\},PC(J)=\{u:J\to R|u\in C(J^{\ast}),u(t_i^+) ,u(t_i^-)\quad exist, u(t_i^-)=u(t_i),i=1,2,\cdots,p\}$.
$PC^1(J)=\{u\in PC(J):|u|_{(t_i,t_{i+1})}\in C^1(t_i,t_{i+1}),u'(t_i^+) ,u'(t_i^-)\quad exist, u'(t_i^-)=u'(t_i),i=1,2,\cdots,p\}$

Lijing Chen and Jitao Sun [31] discussed the nonlinear boundary value problem of first order impulsive functional differential equations by using the upper and lower solution method and monotone iterative technique:
\begin{align*}
	\left\{
	\begin{array}{ll}
		x''(t)=f(t,x(t),x(\theta(t))),& t\in J=[0.T],t\neq t_k,k=1,2,\cdots,p,\\
		\Delta x(t_k)=I_k(x(t_k)),&k=1,2,\cdots,p,\\
		g(x(0),x(T))=0,
	\end{array} 
	\right.
\end{align*}
where $f\in C(J\times R^2,R),I_k\in C(R,R),g\in C(R\times R,R),\theta\in C(J,J).\Delta x(t_k)=x(t_k^+)-x(t_k^-)$, in which $x(t_k^+),x(t_k^-)$ denote the right and left limits of $x(t)$ at $t_k,k=1,2,\cdots,p$, which are fixed such that $0<t_1<t_2<\cdots<t_p<T$.Let $J_0=J\backslash\{t_1,t_2,\cdots,t_p\},$ $\tau=max_k\{t_k-t_{k-1},k=1,2,\cdots,p\}$,here $t_0=0,t_{p+1}=T$.And  $PC(J)=\{u:J\to R|u\in C(J_0),u(t_i^+) ,u(t_i^-)\quad exist, u(t_i^-)=u(t_i),i=1,2,\cdots,p\}$.

The first-order Hamiltonian system has not been involved by many people.
In recent years, variational method has been used by many scholars to study the solutions of differential equation. In fact, it is very difficult to get a strong solution of a differential equation. The general method is to transform the differential equation into an integral equation, and then get its corresponding energy functional. In this way, we can use variational method and critical point theory to study differential equation. Many scholars have done a lot of work on differential equation by using variational method and critical point theory, such as [1, 13, 24, 29]. For the case of differential equations with fixed impulses see [2, 4, 9, 22, 14, 25].For example, Jingli Xie, Jianli Liand Zhiguo Luo gives periodic and subharmonic solutions of second-order Hamiltonian equation with fixed pulse in [1] by using the linking
theorem:
\begin{align*}
	\left\{
	\begin{array}{ll}
		-q''(t)=\nabla F(t,q(t)),& t\neq t_j,t\in R\\
		\Delta q'(t_j)=-g_j(q(t_j)),&j\in Z,
	\end{array} 
	\right.
\end{align*}
where $q\in R^{N},\nabla F(t,q)=grad_qF(t,q),g_j(q)=grad_qG_j(q),G_j\in(R^N,R)$ for each $j\in Z$, and the operator $\Delta$ is defined as $\Delta \dot{q}(t_j)=\dot{q}(t^+_j)-\dot{q}(t^-_j)$, where $\dot{q}(t^
+_j)$ $\dot{q}(t^-_j)$ denotes the right-hand (left-hand) limit of $\dot{q}$ at $t_j$. There exist an $m\in N$ and a $T>0$ such that $0=t_0<t_1<t_2<\cdots<t_m=T, t_{j+m}=t_j+T$, and $g_{j+m}=g_j,j\in N. F:R\times R^{N}\to R$ is $T$-periodic in its first variable and satisfies:$F(t,q)$ is measurable in $t$ for each $q\in R^{N}$ and continuously differentiable in $q$ for a.e. $t\in[0,T]$.

The stochastic pulse differential equation we study is more general than the fixed pulse differential equation above, and it has a wide range of applications, and can simulate the real life situation more. Inspired by [3, 32, 33], we obtain a class of first-order Hamiltonian systems and decide to study the existence of periodic solutions for first-order Hamiltonian systems with impulses under Dirichlet boundary value conditions.\begin{align*}
	u'(t)=A(t)u(t)=JD(t)u(t)=J\nabla H(t,u(t))
\end{align*}

In the motion of Gas Block Simulation of Linear Convection, the circulation generated by the horizontal temperature gradient between the equator and the polar region is accompanied by the rising motion of the equatorial air and the sinking motion of the polar air. Because of the disturbance, the vertical displacement needs to consider the instantaneous impulse.
Therefore, it is reasonable to add such an impulsive condition.

Hence, we consider the existence of solutions to the random pulse linear hamiltonian system boundary value problem:
\begin{equation}\label{1.1}
	\left\{
	\begin{array}{ll}
		u'(t)=A(t)u(t)=JD(t)u(t)=J\nabla H(t,u(t)),&t\in[0,T]\backslash\{\xi_1,\xi_2,\cdots\},\\
		\Delta u(\xi_j)=u(\xi_j^+)-u(\xi_j^-)=b_j(\tau_j),&\xi_j\in(0,T),~j=1,2,\cdots,\\
		u(0)=u(T)=0,
	\end{array} 
	\right.
\end{equation}
where $\forall t\in[0,T],A(t)\in M_{2n}(\mathbb{R})$ is a hamiltonian matrix. $H:[0,T]\times\mathbb{R}^{2n}\to\mathbb{R}$,$(t,u)\to H(t,u)$ is a smooth Hamiltonian.
$\nabla H$ is the gradient of $H$ with respect to $u$.
$D\in gl(2n,\mathbb{R})$ is a symmetric matrix with respect to $t$ continuity,$gl(2n,\mathbb{R})$ is the set of all $2n\times2n$ matrices in the field $\mathbb{R}$.
$u(t):[0,T]\times\Omega\to\mathbb{R}^{2n}$ is a stochastic process. $\forall j\in\mathbb{N},\tau_j:\Omega\rightarrow F_j$, where $F_j:=(0,d_j)$ is a random variable,with $0<d_j<+\infty,$ and $\forall i,j\in\mathbb{N},\tau_i,\tau_j$ are mutually independent when $i\neq j$, $b_j:F_j\to\mathbb{R}^{2n}.$ Set $\xi_{j+1}=\xi_j+\tau_j$. $\{\xi_j\}$ is a strictly increasing random variable sequence i.e. $0<\xi_1<\xi_2<\cdots<\xi_k<\cdots<T$. $u(\xi_j^+):=\lim\limits_{t\to\xi_j^+}u(t),u(\xi_j^-):=\lim\limits_{t\to\xi^-_{j}}u(t)$ under the meaning of the sample orbit. The definition of this defintion are reasonable because $\{\xi_k\}$ will become a series of fixed points under the realization of each sample orbit. We suppose that $\{N(t):t\geq 0\}$ is the simple counting process generated by $\{\xi_k\}$, that is, $\left\{N(t)\geq n\right\}=\left\{\xi_n\leq t\right\}$, and denote $\psi_t$ the $\sigma$-algebra generated by $\left\{N(t),t\geq 0\right\}$.$J$ is the $2n\times2n$ matrices,
\begin{displaymath}
	J=\begin{bmatrix}
		0 & I_n\\
		-I_n & 0
	\end{bmatrix}
\end{displaymath}
where,$I$ is the $n\times n$ identity matrix.

\section{Preliminaries \label {Sec2}}
Let $(\Omega,\psi,P)$ be a probability space. Let $L^q([0,T]\times \Omega,R)$ be the collection of all strongly measurable, $q$th-integrable,$\psi_t$-measurable $R$-valued random variables $x$ with norm $||x||_{L^q}=(E||x||^q)^{\frac{1}{q}}$, where the expectation $E$ is defined by $Ex=\int_{\Omega}xdP$. Let $PC\Big([0,T]\Big):=\Big\{u(t)=u(t,\omega)$ is random process,
$u(\cdot,\omega)$ is a map from $[0,T]$ to $\mathbb{R}^{2n}$ such that $u(t)$ is continuous on $[0,T]\backslash \{\xi_1,\xi_2,\cdots\}$and ${u}(\xi_j^+), {u}(\xi_j^-)$ exist,
$j=1,2,\cdots; u(0)=u(T)=0\Big\}$,$PC\Big([0,T]\Big)$ is a Banach space with norm $\|u\|_{PC}=(\underset{t\in[0,T]}{\max}E|u(t)|^2)^\frac{1}{2}$. Define the Banach space $PC_1 = PC_1\Big([0,T]\Big):=\Big\{u(t)=u(t,\omega)$ is random process,
$u(\cdot,\omega)$ is a map from $[0,T]$ to $\mathbb{R}^{2n}$ such that $u(t)$ and $u^{'}(t)$ is continuous on $[0,T]\backslash \{\xi_1,\xi_2,\cdots\}$ and ${u^{'}}(\xi_j^+), {u^{'}}(\xi_j^-)$ exist,
$j=1,2,\cdots; u(0)=u(T)=0\Big\}$,
with the norm $\|u\|_{PC_1}=\max\Big\{\|u\|_{PC},\|\dot{u}\|_{PC}\Big\}$.\\

\noindent

We introduce the Legendre transformation $H^{\ast}(t,\cdot)$ of $H(t,\cdot)$ and define it as
\begin{align}\label{2.1}
	H^{\ast}(t,v)=(v,u)-H(t,u)
\end{align} 

where
\begin{align*}
	v=\nabla H(t,u),\qquad
	u=\nabla H^{\ast}(t,v)
\end{align*}

\begin{definition}(Fr\'echet derivative)
     $E$ is Banach space, $I:E\to R$ is a functional on $E$, $u\in E$. If there is $A(u)\in E^{\ast}$, such that
     \begin{align*}
     	I(u+\varphi)=I(u)+(A(u),\varphi)+\omega(u,\varphi)
     \end{align*}
 where $(A(u),\varphi)=A(u)\varphi$, represents the value of functional $A(u)$ at $\varphi$, $\omega(u,\varphi)=\circ||\varphi||$, i.e.
 \begin{align*}
 	\lim\limits_{||\varphi||\to 0}\frac{||\omega(u,\varphi)||}{||\varphi||}=0,
 \end{align*}
The functional $I$ is called Fr\'echet differentiable at $u$, $A(u)$ is called the Fr\'echet derivative of $I$ at $u$, then
\begin{align*}
	I(u+\varphi)=I(u)+(I'(u),\varphi)+\circ||\varphi||
\end{align*}
If $I$ is Fr\'echet differentiable for any $u$, denote as $I\in C^1(E,R)$.
\end{definition}

\begin{lemma}
	If $q>2,\alpha>0$ exists and The function $H^{\ast}$ is differentiable in the domain of definition,with $\frac{1}{p}+\frac{1}{q}=1,M^{\ast}=\underset{||v||_{PC}=1}\max H^{\ast}(v),\alpha^{\ast}=(\alpha q)^{-\frac{p}{q}}/p$, make the following conditions hold
	\begin{align*}
		qH(u)&\leq(\nabla H(u),u)\\
		H(u)&\leq\alpha||u||^{q}_{PC}
	\end{align*}
	then,we have
	\begin{align*}
		&pH^{\ast}\geq (\nabla H^{\ast}(v),v)\\
		&H^{\ast}\leq M^{\ast}||v||^{p}_{PC},\quad ||v||_{PC}\geq 1\\
		&H^{\ast}\geq \alpha^{\ast}||v||^{p}_{PC},\\
	\end{align*}
\end{lemma}

\begin{proof}(i)\quad We know that
	\begin{align*}
		v=\nabla H(u)\Rightarrow (\nabla H(u),u)&=(v,u)\geq qH(u)\\
		H^{\ast}(v)=(v,u)-H(u)&\geq (1-\frac{1}{q})(v,u)=\frac{1}{p}(\nabla H^{\ast}(v),v)
	\end{align*}
	(ii)\quad Fixed $v$, let $f(w)=H^{\ast}(wv)$,by condition we have $pf(w)\geq wf^{'}(w)$,if $w\geq 1,w^{p}f(1)\geq f(w)$,i.e.$w^{p}H^{\ast}(v)\geq H^{\ast}(wv)$.When $||v||_{PC}\geq 1$,then $H^{\ast}(\frac{v}{||v||_{PC}})\geq ||v||^{-p}_{PC}H^{\ast}(||v||_{PC}\frac{v}{||v||_{PC}})$.\\
	(iii)\quad By $H^{\ast}(v)=(v,u)-H(u)$ and $H(u)\leq \alpha ||u||^{q}_{PC}$,then
	\begin{align*}
		H^{\ast}(v)&\geq (v,u)-\alpha ||u||^{q}_{PC}\\
		H^{\ast}(v)&\geq \underset{u}\sup((v,u)-\alpha ||u||^{q}_{PC})\\
		&=||v||^{p}_{PC}(\alpha q)^{-\frac{p}{q}}/p.
	\end{align*}
\end{proof}

\begin{lemma}
	If $u$ and $v$ are two random processec,where the expectation $E$ is definded $E(x)=\int_{\Omega}xdP$, then
	\begin{align*}
		(E|u||v|)^2\leq E|u|^2E|v|^2
		\end{align*}
\end{lemma}
\begin{proof}
\begin{align*}
	(E|u||v|)^2&=(\int_{\Omega}|u||v|dP)^2\\
&\leq \int_{\Omega}|u|^2dP\int_{\Omega}|v|^2dP\\
&=E|u|^2E|v|^2
\end{align*}
\end{proof}

\begin{lemma}(Embedding theorem)
	If $u(t)$ is a stochastic process and we denote by $PC_1\to PC$ the embeddinng, then there is a constant $K$ such that
	\begin{align}
		||u||_{PC_1}\leq K||\dot{u}||_{PC}
	\end{align} 
\end{lemma}

\begin{theorem}
	(Mountain path lemma) E is a Banach space, $\varphi\in C^{1}(E,R)$, if $\varphi$ satisfies 
	
	(i)\quad $\varphi(0)=0,{\exists}\rho >0,s.t. \varphi_{\partial B_{\rho}(0)}\geq \alpha >0$;
	
	(ii)\quad $\exists e\in E\backslash \overline{B_{\rho}(0)}, s.t. \varphi(e)\leq 0$.
	
	(iii)\quad the $P.-S.$ condition is fulfilled.
	
	Then, $\varphi$ exists a critical point $u$ satisfying $\varphi^{'}(u)=0$ and $\varphi(u)>max \left\{\varphi(0),\varphi(e)\right\}$.
\end{theorem}

\begin{remark}
	($P.-S.$ condition)\quad Suppose $\varphi\in C^{1}(E,R)$. If $\left\{\varphi(u_k)\right\}$ is bounded and $\left\{\varphi'(u_k)\right\}\to 0$ in $E^{\ast}$, when $k\to\infty$ implies that each $\left\{u_k\right\}$ is sequentially compact set in $E$. Then we call $\varphi$ satisfies $P.-S.$ condition.
\end{remark}

Now, we present some important conclusions that will be used in the next section.
\begin{theorem}
	 Define the function $\forall u\in PC_1$
	\begin{equation}\label{2.8}
		\varphi(u)=E\left[-\dfrac{1}{2}\int_{0}^{T}(J\dot{u},u(t))dt-\int_{0}^{T}H(t,u(t))dt-\dfrac{1}{2}\sum_{k=1}^{\infty}\left(\sum_{j=1}^{k}(Ju(\xi_j),b_j(\tau_j))I_A(\{\xi_j\}_{j=1}^k)\right)\right],
	\end{equation}
	where 
	\begin{displaymath}
		J=\begin{bmatrix}
			0 & I_n\\
			-I_n & 0
		\end{bmatrix},\qquad
		I_{A}(x)=\left\{ \begin{array}{ll}
			1,  & \textrm{if  x~$\in A$ ,}\\
			0,  & \textrm{if  x~$\notin A$ },
		\end{array} \right.\quad (u,v):=\sum_{i=1}^{2n}u(i)v(i)
	\end{displaymath}
	and A is the set consisting of all sample orbits, and $\{\xi_i\}_{i=1}^k$ is a sample orbit.
	
	If $u=Jv$, then
	\begin{align*}\label{2.8}
		\varphi(u) &=E\left[\dfrac{1}{2}\int_{0}^{T}(\dot{v}(t),u(t))dt-\int_{0}^{T}H(t,u(t))dt+\dfrac{1}{2}\sum_{k=1}^{\infty}\left(\sum_{j=1}^{k}(v(\xi_j),b_j(\tau_j))I_A(\{\xi_j\}_{j=1}^k)\right)\right]\\
		&=E\left[-\dfrac{1}{2}\int_{0}^{T}(\dot{v}(t),u(t))dt+\int_{0}^{T}(\dot{v}(t),u(t))dt-\int_{0}^{T}H(t,u(t))dt\right]\\
		&+\left[\dfrac{1}{2}\sum_{k=1}^{\infty}\left(\sum_{j=1}^{k}(v(\xi_j),b_j(\tau_j))I_A(\{\xi_j\}_{j=1}^k)\right)\right]\\
		&=E\left[\dfrac{1}{2}\int_{0}^{T}(J\dot{v}(t),v(t))dt+\int_{0}^{T}(\dot{v}(t),u(t))dt-\int_{0}^{T}H(t,u(t))dt\right]\\
		&+\left[\dfrac{1}{2}\sum_{k=1}^{\infty}\left(\sum_{j=1}^{k}(v(\xi_j),b_j(\tau_j))I_A(\{\xi_j\}_{j=1}^k)\right)\right]
	\end{align*}
	If $H^{\ast}$ in \eqref{2.1} is used to replace $(\dot{v},u)-H(t,u)$, the conjugate action of $T$ periodic function space is obtained.
	\begin{align}
		\chi(v)=E\left[\dfrac{1}{2}\int_{0}^{T}(J\dot{v}(t),v(t))dt+\int_{0}^{T}H^{\ast}(t,\dot{v}(t))dt+\dfrac{1}{2}\sum_{k=1}^{\infty}\left(\sum_{j=1}^{k}(v(\xi_j),b_j(\tau_j))I_A(\{\xi_j\}_{j=1}^k)\right)\right]
	\end{align}
	We can prove that $\chi(v)\in C^{1}(PC_1,\mathbb{R})$ and $\forall v\in PC_1,\forall h \in PC_1$,
	\begin{align}\label{2.9}
		(\chi'(v),h)=&E\left[\dfrac{1}{2}\int_{0}^{T}(J\dot{v}(t),h(t))dt+\int_{0}^{T}(\nabla H^{\ast}(t,\dot{v}(t))-\dfrac{1}{2}Jv(t),\dot{h}(t))dt\notag\right]\\&+\left[\dfrac{1}{2}\sum_{k=1}^{\infty}\left(\sum_{j=1}^{k}(h(\xi_j),b_j(\tau_j))I_A(\{\xi_j\}_{j=1}^k)\right)\right]
	\end{align}
Detailed proof of these will be given in Section 3.
\end{theorem}

\begin{theorem}
	If the random impulsive differential equation \eqref{1.1} has a mild solution $u=\nabla H^{\ast}(t,\dot{v}(t))$, then $v$ is a critical point of $\chi(v)$ i.e. $(\chi'(v),h)=0,\forall h\in PC_1$. And if $v\in PC_1,v$ is a critical point of $\chi(v)$, then $u=\nabla H^{\ast}(t,\dot{v}(t))$ is a mild solution of \eqref{1.1}.
	
	Proof: Suppose $0=t_0<t_1<t_2<\cdots<t_k<t_{k+1}=T$, where $t_1,t_2,\cdots,t_k$ is a sample orbit. $\{t_i\}_{i=1}^k\in A$.
	Let $u\in PC_1$ is a mild solution of (1.1). If $u\in C^1([0,T]\backslash\{t_1,t_2,\cdots,t_k\})\cap PC_1$ then $u$ is the solution of \eqref{1.1} and satisfies
	$$\dot{u}(t)=J\nabla H(t,u(t)).$$
	From hypothesis $u=Jv$ and conjugate symmetry, we have
	\begin{align*}
		J\dot{v}(t)&=J\nabla H(t,u(t))\\
		\dot{v}(t)&=\nabla H(t,u(t))\\
		u(t)&=\nabla H^{\ast}(t,\dot{v}(t)).
	\end{align*}
	Let's take the inner product by $h(t)\in PC_1\cap C^{1}([0,T]\times\Omega,\mathbb{R}^{2n})$ of both sides, then integration from $0$ to $T$, we get
	$$\int_{0}^{T}(u(t),\dot{h}(t))dt=\int_{0}^{T}(\nabla H^{\ast}(t,\dot{v}(t)),\dot{h}(t))dt.$$
	That means
	\begin{equation}\label{2.10}
		\frac{1}{2}\int_{0}^{T}(u(t),\dot{h}(t))dt=\int_{0}^{T}(\nabla H^{\ast}(t,\dot{v}(t))-\frac{1}{2}u(t),\dot{h}(t))dt,
	\end{equation}
	where 
	\begin{align}\label{2.11}
		\int_{0}^{T}(u(t),\dot{h}(t))dt&=\sum_{j=0}^{k}\int_{t_j}^{t_{j+1}}(u(t),\dot{h}(t))dt\nonumber\\
		&=\sum_{j=0}^{k}\nonumber\left[(u(t),h(t))\Big|_{t_j^+}^{t_{j+1}^-}-\int_{t_j}^{t_{j+1}}(\dot{u}(t),h(t))dt\right]\\
		&=-\sum_{j=1}^{k}(h(t_j),\Delta u(t_j))-\int_{0}^{T}(\dot{u}(t),h(t))dt.
	\end{align}
	Then we put \eqref{2.11} into \eqref{2.10} and consider the impulsive condition in \eqref{1.1}, we obtain that
	\begin{equation}
		\frac{1}{2}\int_{0}^{T}(\dot{u}(t),h(t))dt+\int_{0}^{T}(\nabla H^{\ast}(t,\dot{v}(t))-\frac{1}{2}u(t),\dot{h}(t))dt+\frac{1}{2}\sum_{i=1}^{k}(h(t_j),b_j (\tau_j))=0.
	\end{equation}
	By\quad$u(t)=Jv(t)$
	\begin{equation}
		\frac{1}{2}\int_{0}^{T}(J\dot{v}(t),h(t))dt+\int_{0}^{T}(\nabla H^{\ast}(t,\dot{v}(t))-\frac{1}{2}Jv(t),\dot{h}(t))dt+\frac{1}{2}\sum_{j=1}^{k}(h(t_j),b_j (\tau_j))=0.
	\end{equation}
	
	Thus we know $(\chi'(v),h)=0$, i.e. $v$ is a critical point of $\chi(v)$.
	
	On the other hand, if $v\in PC_1$ is a critical point of $\chi$, i.e. $(\chi'(v),h)=0, \forall h\in PC_1,$
	\begin{align}
			&\frac{1}{2}\int_{0}^{T}(J\dot{v}(t),h(t))dt+\int_{0}^{T}(\nabla H^{\ast}(t,\dot{v}(t))-\frac{1}{2}Jv(t),\dot{h}(t))dt+\frac{1}{2}\sum_{j=1}^{k}(h(t_j),b_j (\tau_j))=0
		\end{align}
	by $u=Jv$
	\begin{align}
		&\frac{1}{2}\int_{0}^{T}(\dot{u}(t),h(t))dt+\int_{0}^{T}(\nabla H^{\ast}(t,\dot{v}(t))-\frac{1}{2}u(t),\dot{h}(t))dt+\frac{1}{2}\sum_{j=1}^{k}(h(t_j),b_j (\tau_j))=0, \forall h\in S.
	\end{align}
Since $h\in C^1$, we know $h(t^+_j)=h(t^-_j), j=1,2,3,\cdots,$ and $h(0)=h(T)=0$. We will prove $u=Jv$ is the solution of \eqref{1.1}:
	\begin{align*}
		\dfrac{1}{2}\sum_{i=0}^{k}\int_{t_i}^{t_{i+1}}(\dot{u}(t),h(t))dt+\int_{0}^{T}(\nabla H^{\ast}(t,\dot{v}(t))-\frac{1}{2}u(t),\dot{h}(t))dt+\frac{1}{2}\sum_{i=1}^{k}(h(t_j),b_j (\tau_j))=0.
	\end{align*}
For the convenience, let $t_0=0,t_{k+1}=T$ and $h(t_0)=h(t_{k+1})=0$,
	\begin{align*}
		\dfrac{1}{2}\sum_{i=0}^{k}\nonumber\int_{t_i}^{t_{i+1}}(\dot{u}(t),h(t))dt+\int_{0}^{T}(\nabla H^{\ast}(t,\dot{v}(t))-\frac{1}{2}u(t),\dot{h}(t))dt+\frac{1}{2}\sum_{i=1}^{k}(h(t_j),b_j (\tau_j))=0.
	\end{align*}
	\begin{align}\label{2.13}
		\Rightarrow\int_{0}^{T}(\dot{h},u-\nabla H^{\ast}(t,\dot{v}(t)))dt+\dfrac{1}{2}\sum_{j=1}^{k}(h(t_j),\Delta u(t_j)-b_j(\tau_j))=0
	\end{align}
	Set
	\begin{displaymath}
		\delta_j(t)=\left\{ \begin{array}{ll}
			1,  & \textrm{if  t~$=t_j$ ,}\\
			0,  & \textrm{if  t~$\ne t_j$ },
		\end{array} \right.
	\end{displaymath}
	then $\eqref{2.13}$ can be written as
	\begin{align}\label{2.2}
		\int_{0}^{T}\left[(\dot{h},u-\nabla H^{\ast}(t,\dot{v}(t)))+\dfrac{1}{2}\sum_{j=1}^{k}(h(t_j),\Delta u(t_j)-b_j(\tau_j))\delta_j(t)\right]dt=0
	\end{align}
	thus $u\in PC_1$ is the mild solution of the equation
	\begin{align*}
		u&=\nabla H^{\ast}(t,\dot{v}(t)),\quad t\in [0,T]\backslash\{t_1,t_2,\cdots,t_k\}\\
			\dot{u}(t)&=J\nabla H(t,u(t)),\quad t\in [0,T]\backslash\{t_1,t_2,\cdots,t_k\}
	\end{align*}
	and \eqref{2.2} imply that the random impulsive condition $\Delta u(t_j)=b_j(\tau_j), j=1,2,\cdots,k$ hold.\\
	Thus $u\in PC_1$ is a mild solution of \eqref{1.1}.

\end{theorem}

\section{Main Results \label {Sec3}}

\begin{theorem}\label{thm3.1}
	When $b_j(\tau_j)$ satisfy the followng assumptions respectively
	
	(H1) Let $B=E\Big(\sum_{k=1}^\infty(\sum_{j=1}^k|b_j(\tau_j)|I_A(\{\xi_j\}_{j=1}^k))\Big)<+\infty$.
	
	Then the $\chi(v)$ defined in \eqref{2.8} fulfills $\chi\in C^1(PC_1,\mathbb{R})$ and satisfies \eqref{2.9}.
\end{theorem}

\begin{proof}
	We divide the proof into several parts.
	1. Let $J_1(v)=\dfrac{1}{2}\int_{0}^{T}E(J\dot{v}(t),v(t))dt$, and we will prove that $J_1(v)\in C^1(PC_1,\mathbb{R}).$
	
	$\forall v,h\in PC_1$,we have 
	$$J_1(v+h)=J_1(v)+\dfrac{1}{2}\int_{0}^{T}E(J\dot{v},h)dt+\dfrac{1}{2}\int_{0}^{T}E(J\dot{h},v)dt+\dfrac{1}{2}\int_{0}^{T}E(J\dot{h},h)dt.$$
	Since \begin{align*}
		\left(\dfrac{1}{2}\int_{0}^{T}E(J\dot{h},h)dt\right)^2
		&\leq(\dfrac{1}{2}\int_{0}^{T}E|J\dot{h}||h|dt)^2\\
		&\leq T\dfrac{1}{4}\int_{0}^{T}(E|\dot{h}||h|)^2dt\\
		&\leq T\dfrac{1}{4}\int_{0}^{T}E|\dot{h}|^2E|h|^2dt\\
		&\leq T\dfrac{1}{4}\int_{0}^{T}\underset{t\in[0,T]}{\max}E|\dot{h}|^2\underset{t\in[0,T]}{\max}E|h|^2dt\\
		&=T^2\frac{1}{4}||h||^2_{PC}||\dot{h}||^2_{PC}\\
		&\leq T^2\frac{1}{4}||h||^4_{PC_1}
	\end{align*}
	then we have \begin{align*}
		\left|\dfrac{1}{2}\int_{0}^{T}E(J\dot{h},h)dt\right|\leq T\frac{1}{2}||h||^2_{PC_1}\\
		\lim\limits_{\|h\|_{PC_1}\to 0}\dfrac{\frac{1}{2}\int_{0}^{T}E(J\dot{h},h)dt}{\|h\|_{PC_1}}=0.
	\end{align*}
	It follows that $$(J'_1(v),h)=\dfrac{1}{2}\int_{0}^{T}E(J\dot{v},h)dt+\dfrac{1}{2}\int_{0}^{T}E(J\dot{h},v)dt.$$
	For fixed $v$,$J'_1(v)$ is a linear functional with respect to $h$. By Cauchy-Schwarz inequality, we have 
	\begin{align*}
		(\dfrac{1}{2}\int_{0}^{T}E(J\dot{v},h)dt)^2
		&\leq(\dfrac{1}{2}\int_{0}^{T}E|J\dot{v}||h|dt)^2\\
		&=T\dfrac{1}{4}\int_{0}^{T}(E|\dot{v}||h|)^2dt\\
		&\leq T\dfrac{1}{4}\int_{0}^{T}E|\dot{v}|^2E|h|^2dt\\
		&\leq T\dfrac{1}{4}\int_{0}^{T}\underset{t\in[0,T]}{\max}E|\dot{v}|^2\underset{t\in[0,T]}{\max}E|h|^2dt\\
		&=T^2\frac{1}{4}||h||^2_{PC}||\dot{v}||^2_{PC}\\
		&\leq T^2\frac{1}{4}||h||^2_{PC_1}||v||^2_{PC_1}
	\end{align*}
	then, we have $(\dfrac{1}{2}\int_{0}^{T}E(J\dot{h},v)dt)^2
	\leq T^2\frac{1}{4}||h||^2_{PC_1}||v||^2_{PC_1}$
	\begin{align*}
		|(J'_1(v),h)|
		&\leq\dfrac{T}{2}\|h\|_{PC_1}\|v\|_{PC_1}+\dfrac{T}{2}\|v\|_{PC_1}\|h\|_{PC_1}\\
		&=(T\|h\|_{PC_1})\|v\|_{PC_1},
	\end{align*}
	where $T\|v\|_{PC_1}$ is independent of $h$, therefore $J'_1(v)$ is a bounded functional in $PC_1$.
	
	2. Let $J_2(v)=\int_{0}^{T}EH^{\ast}(t,\dot{v}(t))dt$, and we will prove that $J_2(v)\in C^1(PC_1,\mathbb{R}).$
	
	$\forall v,h\in PC_1$, we have \begin{align*}
		J_2(v+h)&=J_2(v)+\int_{0}^{T}E(u,\dot{h}(t))dt\\
		&=J_2(v)+\int_{0}^{T}E(\nabla H^{\ast}(t,\dot{v}(t)),\dot{h}(t))dt
	\end{align*}
	
	therefore when $v$ fixed, $J'_2(v)$ is a inear functional w.r.t. $h$.
	\begin{align*}
		(J'_2(v),h)&=\int_{0}^{T}E(\nabla H^{\ast}(t,\dot{v}(t)),\dot{h}(t))dt\\
		&\leq(\int_{0}^{T}E\mid\nabla H^{\ast}(t,\dot{v}(t))\mid^2 E\mid\dot{h}(t)\mid^2 dt)^{\frac{1}{2}}
	\end{align*}
	By $u=Jv=\nabla H^{\ast}(t,\dot{v}(t))$ ,bring in the above formula,we have
	\begin{align*}
		(J'_2(v),h)&\leq(\int_{0}^{T}E(\mid Jv\mid^2) E\mid\dot{h}(t)\mid^2 dt)^{\frac{1}{2}}\\
		&\leq T||h||_{PC_1}||Jv||_{PC_1}
	\end{align*}

	where $\|v\|_{PC_1}$ is independent of $h$, therefore $J'_2(v)$ is a bounded functional in $PC_1$.
	
	3. Let $J_3(v)=\sum_{k=1}^{\infty}\sum_{j=1}^{k}E(v(\xi_j),b_j(\tau_j))I_A(\{\xi_j\}_{j=1}^k)$, and we will prove that $J_3(v)\in C^1(PC_1,\mathbb{R}).$
	
	$\forall v,h\in PC_1$, we have \begin{align*}
		J_3(v+h)=J_3(v)+\sum_{k=1}^{\infty}\sum_{j=1}^{k}E(h(\xi_j),b_j(\tau_j))I_A(\{\xi_j\}_{j=1}^k),
	\end{align*}
	then $$(J'_3(v),h)=\sum_{k=1}^{\infty}\sum_{j=1}^{k}E(h(\xi_j),b_j(\tau_j))I_A(\{\xi_j\}_{j=1}^k).$$
	when $v$ fixed, $J'_2(v)$ is a inear functional w.r.t. $h$. By the property of hamiltonian matrix, we have \begin{align*}
		|(J'_3(v),h)|\leq (\sum_{k=1}^{\infty}\sum_{j=1}^{k}E|h(\xi_j)|^2E|b_j(\tau_j)|^2I_A(\{\xi_j\}_{j=1}^k))^\frac{1}{2}\leq B||h||_{PC_1},
	\end{align*}
	thus $J'_3(v)$ is a bounded functional in $PC_1$.
	
	From parts 1.$\sim$3. we can conclude that $\chi\in C^1(PC_1,\mathbb{R})$ and satisfies \eqref{2.9}.
\end{proof}

\begin{theorem}\label{thm3.2}
	When impulse satisfies assumption (H1), $H(u)$ and $\nabla H(u)$ satisfy the followng assumptions respectively
	
	(H2) $H(u)\leq\alpha\|u\|^{q}_{PC}, (q>2,\frac{1}{p}+\frac{1}{q}=1)$ holds, where $\alpha$ is a constant.
	
	(H3)$(\nabla H(u),u)\geq qH(u)$
	
	Then $\chi(v)$ satisfies $P.-S.$ condition.
\end{theorem}
	
\begin{proof}
	Step 1:We will prove that $\left\{v_k\right\}$ is a bounded  set in $PC_1$ provided $\left\{\chi(v_k)\right\}$ ia a bounded set and $\left\{\chi^{'}(v_k)\right\}\rightarrow0$ in $PC^{\ast}_1$ as $k\rightarrow\infty$.
	
	Since
	\begin{align*}
		\chi(v)&=\dfrac{1}{2}\int_{0}^{T}E(J\dot{v}(t),v(t))dt+\int_{0}^{T}EH^{\ast}(t,\dot{v}(t))dt+\dfrac{1}{2}\sum_{k=1}^{\infty}\left(\sum_{j=1}^{k}E(v(\xi_j),b_j(\tau_j))I_A(\{\xi_j\}_{j=1}^k)\right)\\
		&=-\dfrac{1}{2}\int_{0}^{T}E(\dot{v}(t),Jv(t))dt+\int_{0}^{T}EH^{\ast}(t,\dot{v}(t))dt+\dfrac{1}{2}\sum_{k=1}^{\infty}\left(\sum_{j=1}^{k}E(v(\xi_j),b_j(\tau_j))I_A(\{\xi_j\}_{j=1}^k)\right)\\
	\end{align*}
	by
	\begin{align*}
		Jv=u=\nabla H^{\ast}(t,\dot{v}(t))
	\end{align*}
	then
	\begin{align*}
		\chi(v)&=-\dfrac{1}{2}\int_{0}^{T}E(\dot{v}(t),\nabla H^{\ast}(t,\dot{v}(t))dt+\int_{0}^{T}EH^{\ast}(t,\dot{v}(t))dt\\&+\dfrac{1}{2}\sum_{k=1}^{\infty}\left(\sum_{j=1}^{k}E(v(\xi_j),b_j(\tau_j))I_A(\{\xi_j\}_{j=1}^k)\right)\\
		&\geq (1-\frac{p}{2})\int_{0}^{T}EH^{\ast}(t,\dot{v}(t))dt+\dfrac{1}{2}\sum_{k=1}^{\infty}\left(\sum_{j=1}^{k}E(v(\xi_j),b_j(\tau_j))I_A(\{\xi_j\}_{j=1}^k)\right)\\
		&\geq (1-\frac{p}{2})\alpha^{\ast}T||\dot{v}||^{p}_{PC}-\frac{1}{2}B||v||_{PC_1}
	\end{align*}
	and $K^{p}||\dot{v}||^{p}_{PC}\geq ||v||^{p}_{PC_1}$, then
	\begin{align*}
		\chi(v)\geq (1-\frac{p}{2})\alpha^{\ast}\frac{T}{K^{p}}||v||^{p}_{PC_1}-\frac{1}{2}B||v||_{PC_1},1<p<2.
	\end{align*}
	
	Then put ${v_k}$ into the above formula, and due to boundedness of ${\chi(v_k)}$, we have
	\begin{align*}
		+\infty>|\chi(v_k)|\geq (1-\frac{p}{2})\alpha^{\ast}\frac{T}{K^{p}}||v_k||^{p}_{PC_1}-\frac{1}{2}B||v_k||_{PC_1},1<p<2.
	\end{align*}
	If ${v_k}$ is unbounded, then there is subsequence $\lbrace{v_{k_i}}\rbrace\subset\lbrace{v_k}\rbrace$, such that $||v_{k_i}||_{PC_1}\to+\infty$, so then the right end of the above formula tends to infinity, which is contradict with the boundedness of $\lbrace\chi({v_k})\rbrace$. 
	
		Step2:we show that $\left\{v_k\right\}$ is a bounded sequence in $PC_1$,and $\chi^{'}(v_k)\rightarrow0,(k\rightarrow\infty)$ in $PC_1^{\ast}$,then $\left\{v_k\right\}$ is a sequetial compact set in $PC_1$.

We consider that, $H(t,u)$ is a smooth Hamiltonian and due to $\left\{v_k\right\}$ is a bounded set in $PC_1$, then there is $\left\{v_{k_i}\right\}\subset\left\{v_k\right\}$, satisfies

$v_{k_i}\to v$ isuniform convergence in $\left(0,T\right]$;

$v_{k_i}\rightharpoonup v$ is weak converge in $PC_1$.
\begin{align*}
	\Rightarrow 
	||v_{k_i}-v_{k_j}||_{PC_1}\leq K||\dot{v}_{k_i}-\dot{v}_{k_j}||_{PC}=||\nabla H(Jv_{k_i}-Jv_{k_j})||_{PC}\to 0,\quad (i,j\to\infty)
\end{align*}

From this, when $i,j\to \infty,\Rightarrow ||v_{k_i}-v_{k_j}||_{PC_1}\to 0$, we have $\left\{v_{k_j}\right\}$ is a Cauchy sequential compact in $PC_1$, and by the completeness
of $PC_1$, we further get {uki} $\left\{v_{k_j}\right\}$ convergent in $PC_1$. Then $\left\{v_{k}\right\}$ is a sequential compact in $PC_1$.

	From above, we can deduce $\chi$ satisfies $P.-S.$ condition on $PC_1$.
\end{proof}

\begin{theorem}\label{thm3.3} Suppose that $H(u)$ and $\nabla H^{\ast}(u)$ satisfies:
	
	(H1)\quad $B=E\Big(\sum_{k=1}^\infty\sum_{j=1}^k|b_j(\tau_j)|I_A(\{\xi_j\}_{j=1}^k)\Big)<+\infty$
	
	(H2)\quad $qH(u)\leq (\nabla H(u),u)$, $\forall(t,u)\in[0,T]\times R^{2n}$.

	(H3)\quad $H(u)\leq \hat{\alpha} ||u||^{q}_{PC},q>2,\hat{\alpha}>0,\frac{1}{p}+\frac{1}{q}=1,\hat{\alpha}^{\ast}=(\hat{\alpha} q)^{-\frac{p}{q}}/p,(1-\frac{p}{2})\hat{\alpha}^{\ast}>\frac{1}{2}B$.
	
	Then $\chi(v)\in(PC_1,\mathbb{R})$ and $\chi(v)$ satisfies $P.-S.$ condition in $PC_1$. By Mountain pass lemma, we can get $\chi$ has a critical point, i.e. Equation \eqref{1.1} has at least a mild solution in $PC_1$.

\end{theorem}

\begin{proof}
	(1)\quad By hypothesis (H1), and using similar approach with Theoerm 3.1 ,we can prove that $\chi(v)\in C^{1}(PC_1,R)$.\\
	(2) Next we prove $\chi(v)$ satisfies $P.-S.$ condition on $E$.\\
	1) We will prove if $\left\{\chi(v_k)\right\}$ is a bounded set and $\chi^{'}(v_k)\to 0,(k\to \infty)$ in $PC_1^{\ast}$, then ${v_k}$ is a bounded set in $PC_1$.
	\begin{align*}
		\chi(v)&=\dfrac{1}{2}\int_{0}^{T}E(J\dot{v}(t),v(t))dt+\int_{0}^{T}EH^{\ast}(t,\dot{v}(t))dt+\dfrac{1}{2}\sum_{k=1}^{\infty}\left(\sum_{j=1}^{k}E(v(\xi_j),b_j(\tau_j))I_A(\{\xi_j\}_{j=1}^k)\right)
	\end{align*}
	\begin{align*}	
		(\chi'(v),h)&=\dfrac{1}{2}\int_{0}^{T}E(J\dot{v}(t),h(t))dt+\int_{0}^{T}E(\nabla H^{\ast}(t,\dot{v}(t))-\dfrac{1}{2}Jv(t),\dot{h}(t))dt\\
		&+\dfrac{1}{2}\sum_{k=1}^{\infty}\left(\sum_{j=1}^{k}E(h(\xi_j),b_j(\tau_j))I_A(\{\xi_j\}_{j=1}^k)\right)
	\end{align*}
	then
	\begin{align*}
		\chi(v_k)-(\chi'(v_k),v_k)&=\int_{0}^{T}EH^{\ast}(t,\dot{v}_{k}(t))dt-\int_{0}^{T}E(\nabla H^{\ast}(t,\dot{v}_{k}(t))-\dfrac{1}{2}Jv_{k}(t),\dot{v}_{k}(t))dt\\
		&=\int_{0}^{T}EH^{\ast}(t,\dot{v}_{k}(t))dt-\int_{0}^{T}E(\frac{1}{2}\nabla H^{\ast}(t,\dot{v}_{k}(t)),\dot{v}_{k}(t))dt\\
		&\geq \int_{0}^{T}EH^{\ast}(t,\dot{v}_{k}(t))dt-\frac{p}{2}\int_{0}^{T}EH^{\ast}(t,\dot{v}_{k}(t))dt\\
		&= (1-\frac{p}{2})\int_{0}^{T}EH^{\ast}(t,\dot{v}_{k}(t))dt\\
		&\geq (1-\frac{p}{2})\hat{\alpha}^{\ast}||\dot{v_k}||^{p}_{PC}
	\end{align*}
	then if $\left\{v_k\right\}$ is unbounded, by$||v||_{PC_1}\leq K||\dot{v}||_{PC}$, then $\left\{\dot{v}_k\right\}$ is unbounded,then it's in contradiction with $\chi(v_k)$ is bounded with $\chi^{'}(v_k)\to 0,(k\to \infty)$.
	
	2)We will prove $\left\{v_k\right\}$ is a bounded set in $PC_1$, $\chi^{'}(v_k)\to 0(k\to \infty)$ in $PC_1^{\ast}$, then $\left\{v_k\right\}$ is a sequential compact set in $PC_1$.
	
	We consider that, $H(t,u)$ is a smooth Hamiltonian and due to $\left\{v_k\right\}$ is a bounded set in $PC_1$, then there is $\left\{v_{k_i}\right\}\subset\left\{v_k\right\}$, satisfies
	
	$v_{k_i}\to v$ is uniform convergence in $\left(0,T\right]$;
	
	$v_{k_i}\rightharpoonup v$ is weak converge in $PC_1$.
	\begin{align*}
		\Rightarrow 
		||v_{k_i}-v_{k_j}||_{PC_1}\leq K||\dot{v}_{k_i}-\dot{v}_{k_j}||_{PC}=||\nabla H(Jv_{k_i}-Jv_{k_j})||_{PC}\to 0
	\end{align*}
	
	From this, when $i,j\to \infty,\Rightarrow ||v_{k_i}-v_{k_j}||_{PC_1}\to 0$, we have $\left\{v_{k_j}\right\}$ is a Cauchy sequential compact in $PC_1$, and by the completeness
	of $PC_1$, we further get {uki} $\left\{v_{k_j}\right\}$ convergent in $PC_1$. Then $\left\{v_{k}\right\}$ is a sequential compact in $PC_1$.
	
	By 1),2) we have $\chi(v)$ satisfies $P.-S.$ condition on $PC_1$.
	
	(3)\quad At last, we verify whether $\chi(v)$ fuills the conditions of Mountain path lemma.
	
	a)\quad It is obvious that $\chi(0)=0$.
	Since
\begin{align*}
	\chi(v)&=\dfrac{1}{2}\int_{0}^{T}E(J\dot{v}(t),v(t))dt+\int_{0}^{T}EH^{\ast}(t,\dot{v}(t))dt+\dfrac{1}{2}\sum_{k=1}^{\infty}\left(\sum_{j=1}^{k}E(v(\xi_j),b_j(\tau_j))I_A(\{\xi_j\}_{j=1}^k)\right)\\
	&=-\dfrac{1}{2}\int_{0}^{T}E(\dot{v}(t),Jv(t))dt+\int_{0}^{T}EH^{\ast}(t,\dot{v}(t))dt+\dfrac{1}{2}\sum_{k=1}^{\infty}\left(\sum_{j=1}^{k}E(v(\xi_j),b_j(\tau_j))I_A(\{\xi_j\}_{j=1}^k)\right)\\
\end{align*}
by
\begin{align*}
	Jv=u=\nabla H^{\ast}(t,\dot{v}(t))
\end{align*}
then
\begin{align*}
	\chi(v)&=-\dfrac{1}{2}\int_{0}^{T}E(\dot{v}(t),\nabla H^{\ast}(t,\dot{v}(t))dt+\int_{0}^{T}EH^{\ast}(t,\dot{v}(t))dt\\&+\dfrac{1}{2}\sum_{k=1}^{\infty}\left(\sum_{j=1}^{k}E(v(\xi_j),b_j(\tau_j))I_A(\{\xi_j\}_{j=1}^k)\right)\\
	&\geq (1-\frac{p}{2})\int_{0}^{T}EH^{\ast}(t,\dot{v}(t))dt+\dfrac{1}{2}\sum_{k=1}^{\infty}\left(\sum_{j=1}^{k}E(v(\xi_j),b_j(\tau_j))I_A(\{\xi_j\}_{j=1}^k)\right)\\
	&\geq (1-\frac{p}{2})\hat{\alpha}^{\ast}T||\dot{v}||^{p}_{PC}-\frac{1}{2}B||v||_{PC_1}
\end{align*}
and $K^{p}||\dot{v}||^{p}_{PC}\geq ||v||^{p}_{PC_1}$, then
\begin{align*}
	\chi(v)\geq (1-\frac{p}{2})\hat{\alpha}^{\ast}\frac{T}{K^{p}}||v||^{p}_{PC_1}-\frac{1}{2}B||v||_{PC_1},1<p<2.
\end{align*}
	
	Take $\rho=(\frac{K^P}{T})^{\frac{1}{p-1}}>0,u\in\partial B_{\rho}(0)$, then $||v||_{PC_1}=\rho$, then
	\begin{align*}
			\chi(v)&\geq (1-\frac{p}{2})\hat{\alpha}^{\ast}\frac{T}{K^{p}}||v||^{p}_{PC_1}-\frac{1}{2}B||v||_{PC_1}\\
			&\geq ||v||_{PC_1}((1-\frac{p}{2})\hat{\alpha}^{\ast}\frac{T}{K^{p}}||v||^{p-1}_{PC_1}-\frac{1}{2}B)\\
		&= ||v||_{PC_1}((1-\frac{p}{2})\hat{\alpha}^{\ast}-\frac{1}{2}B)>0
	\end{align*}
	Here we used the assumption (H3).
	
	b)\quad Now, we consider $v_{1}(t)=(\cos\frac{2\pi t}{T})e+(\sin\frac{2\pi t}{T})Je$, where$||e||_{PC}=||e||_{PC_1}\geq \max (1,\rho)$, then
	\begin{align*}
		\int_{0}^{T}(J\dot{v}_{1}(t),v_{1}(t))=-\left(\frac{T^{2}}{2\pi}\right)||e||^{2}_{PC}=-\left(\frac{T^{2}}{2\pi}\right)||e||^{2}_{PC_1}
	\end{align*}
	\begin{align*}
		\chi(v_{1})\leq -\left(\frac{T^{2}}{4\pi}\right)||e||^{2}_{PC_1}+TM^{\ast}||e||^{p}_{PC_1}+\frac{1}{2}B||e||_{PC_1}
	\end{align*}
	since $p<2$ and when $||e||_{PC_1}\to +\infty$ we have $\chi(v_{1})\to -\infty$, there is always an $e$ such that $\chi(v_{1})<0$, by $||e||_{PC_1}>\rho$, so $v_1\in E\setminus\overline{B_{\rho}(0)}$.
	
	From a),b) and $\chi(v)\in C^{1}(PC_1,\mathbb{R})$ satisfies $P.-S.$ condition, by Mountain path lemma, we can obtain $\chi$ has a critical point, i.e. equation \eqref{1.1} has weak solutions in $PC_1$. This completes the proof.

\end{proof}

\section{Example \label {Sec4}}

The main result could have many applications, now, we give an example to illustrate this theorem.
We consider the following  random impulsive differential equation with boundary value
problems.

\begin{equation}
	\label{eq1}
	\left\{
	\begin{array}{ll}
		u'(t)=10|u|^8Ju(t),&t\in[0,T]\backslash\{\xi_1,\xi_2,\cdots\},\\
		\Delta u(\xi_j)=u(\xi_j^+)-u(\xi_j^-)=\frac{j}{4^j}\tau_j,&\xi_j\in(0,T),~j=1,2,\cdots,\\
		u(0)=u(T)=0,
	\end{array} 
	\right.
\end{equation}

Its Hamilton quantity is $H(u)=|u|^{10}$. By transforming it into a polar equation we know that (4.1) period is $\frac{\pi}{5}|u|^{-8}$.

Let $\tau_j\sim U(0,\frac{1}{2^j})$, then the probability density function of $\tau_j$ is 
\begin{equation}
	\label{eq1}
	p(x)=\left\{
	\begin{array}{ll}
		2^j&x\in(0,\frac{1}{2^j}),\\
	 0&x\notin(0,\frac{1}{2^j}),\\
	\end{array} 
	\right.
\end{equation}

Set $\xi_0=0,\xi_{j+1}=\xi_{j}+\tau_{j}$. Obviously, $\left\{\xi_j\right\}$ is a process with independent increments and the
impulsive moments $\xi_j$ form a strictly increasing sequence. And for every $j\in\mathbb{N}$,
\begin{align}
	\xi_j<\xi_{j+1}<\frac{1}{2}+\frac{1}{2^2}+\cdots+\frac{1}{2^{j+1}}<1
\end{align}
So in this example, $b_j(\tau_j)=\frac{j}{4^j}\tau_{j}$, and $\tau_{j}$ is a random variable defined from $\Omega$ to $F_j= (0,d_j) = (0,\frac{1}{2^j}).$
Suppose $\tau_j$ and $\tau_i$ are independent of each other when $i\neq j$, $u(\xi_j^+):=\lim\limits_{t\to\xi_j^+}u(t),u(\xi_j^-):=\lim\limits_{t\to\xi^-_{j}}u(t)$.
\begin{align*}
	\sum_{j=1}^k|b_j(\tau_j)|=\sum_{j=1}^k\frac{j}{4^j}(\tau_j)\leq\frac{4^{-k}}{9}(-3k+4^{k+1}-4)
\end{align*}
So, we have proved that $B=E\Big(\sum_{k=1}^\infty\sum_{j=1}^k|b_j(\tau_j)|I_A(\{\xi_j\}_{j=1}^k)\Big)<\frac{4}{9}$.

\begin{align*}
	(\nabla u,u)&=(10|u|^8u,u)={10}|u|^{10}\geq 10|u|^{10}\\
	H(u)&=|u|^{10}\leq||u||_{PC}^{10}
\end{align*}
then we have $q=10,p=\frac{10}{9},\hat{\alpha}=1$, then $\hat{\alpha}^{\ast}\approx0.696$, then $(1-\frac{p}{2})\hat{\alpha}^{\ast}\geq\frac{1}{2}B$.

So, the equation (4.1) meets all the conditions of the theorem (3.3). By Mountain pass lemma, we can get equation (4.1) has at least a mild solution in $PC_1$.

\section{Conclusion \label {Sec5}}
\quad\quad In this paper, we study the existence of solutions for first-order impulsive Hamiltonian systems (1.1) by using the critical point theory and the variational method. First, we prove that the conjugate action (2.4) of the functional is Fr\'echet differentiable. Then we prove that the critical point of (2.4) is a mild solution of Hamiltonian system (1.1). Then, it is proved that the conjugate action of the functional (2.4) satisfies the P. - S. condition under the given conditions .Finally, through mountain pass lemma, it is found that the critical point of (2.4) is a mild solution of Hamiltonian system (1.1).

In our future research, we will continue to study Hamiltonian systems with random impulses. For example, we will study some properties of second-order Hamiltonian systems with impulses.

\end{document}